\documentclass[12pt]{article}

\usepackage{graphicx}
\usepackage{subfigure}
\usepackage{algorithmic}
\usepackage{algorithm}

\usepackage{amsthm,amsmath}
\newtheorem{thm}{Theorem}
\numberwithin{thm}{section}

\newtheorem{lem}[thm]{Lemma}
\newtheorem{defn}[thm]{Definition}

\usepackage{amssymb,latexsym,amsmath,epsfig}
\usepackage{graphicx,relsize}
\usepackage[numbers,square]{natbib}
\usepackage{mathtools}
\usepackage{mathrsfs}               

\newcommand{\mise}{mis\`{e}re }

\newcommand{\R}{\mathcal{R}}
\renewcommand{\L}{\mathcal{L}}
\newcommand{\N}{\mathcal{N}}
\renewcommand{\P}{\mathcal{P}}

\renewcommand{\v}{\phantom{'} \vert \phantom{'}}

\newcommand{\lopt}{\boldsymbol{L}}
\newcommand{\ropt}{\boldsymbol{R}}

\usepackage{array,color}
\usepackage{colortbl}
\usepackage{tabularx}
\definecolor{gr}{gray}{0.6}

\graphicspath{{Images/}}

\newcommand{\Gb}{\overline{G}}

\usepackage{tikz}

\makeatletter
\def\imod#1{\allowbreak\mkern10mu({\operator@font mod}\,\,#1)}
\makeatother

\begin{document}

\begin{center}
{\bf Restricted developments in partizan mis\`ere game theory}
\ \\ \ \\
REBECCA MILLEY 
\\ {\small {\em Computational Mathematics, Grenfell Campus, Memorial University, Canada.}\\{\em rmilley@grenfell.mun.ca  }}\\
\ \\ GABRIEL RENAULT\footnote{Supported by the ANR-14-CE25-0006 project of the French National Research Agency.}
\\ {\small {\em gabriel.renault@ens-lyon.org}}

\end{center}

\vskip 30pt
 
\begin{abstract}
Mis\`ere games have excited new interest over the past decade, with the introduction of an indistinguishability relation for analyzing positions modulo  restricted subsets of games. 
We present a survey of recent progress in the theory of partizan mis\`ere games, including some results for general mis\`ere play, but focussing primarily on this {\em restricted} mis\`ere play. We discuss new and current work on game comparison and game inverses, as well as  ongoing research around reversibility and canonical forms in restricted mis\`ere play. We also show how general results in each of these areas have been applied to specific games to find solutions under mis\`ere play.
\end{abstract}

\section{Introduction}
Most research in combinatorial game theory  assumes normal play, where the first player unable to move loses, as opposed to mis\`ere play, where the first player unable to move wins.
It is rather remarkable how much changes when we simply   switch the goal from {\em getting} the last move to {\em avoiding} the last move.  At first glance one might think mis\`ere play is merely  the `opposite' of normal play, but this is not at all the case. There is actually no relationship between normal outcome and mis\`ere outcome: for every pair of (not necessarily distinct) outcomes 
$\mathscr{O}_1,\mathscr{O}_2 \in \{\mathscr{L, N, P, R}\}$, there is a game with normal outcome $\mathscr{O}_1$ and mis\`ere outcome $\mathscr{O}_2$ \cite{MesdaO2007}. Likewise, strategies from normal play are in general neither the same nor reversed for mis\`ere play. For example, in normal play, Left would always choose a move to $1 = \{0 | \cdot\}$ over a move to $0=\{\cdot | \cdot\}$, but in mis\`ere play there are situations\footnote{Left wins playing first on the single game $0$ and loses playing first on the single game $1$, but loses playing first on $0+*$ and wins playing first on $1+*$. } in which Left should choose $1$ over $0$ and others where Left prefers $0$ over $1$. This means that $0$ and $1$ are incomparable in mis\`ere play, which goes against our intuition that Left is trying to run out of moves before Right.

So we really are in a fog in mis\`ere play. We look to the elegant algebra of normal-play games and hope for some semblance of structure, but we are dismayed at every turn:

\begin{itemize}

\item Zero is trivial. In normal play we have the wonderful property that every previous-win game is equal to zero. In mis\`ere, the zero game is next-win, but our hopes that perhaps every next-win game is equal to zero are more than dashed: in fact, only the game $\{\cdot | \cdot\}$ is equal to zero \cite{MesdaO2007}. In particular, for any game $G\not = 0$, the game $G-G$ is not equal to zero (a very troublesome fact indeed).  Consequently, there are no non-zero inverses, and there is no longer an easy test for the equality and inequality of games.

\item Equality (and inequality) is rare and difficult to prove. Partly due to the triviality of zero,  equivalence classes induced by the equality relation are much smaller in mis\`ere play, and it is not often possible to simplify games. 
Inequality is  likewise uncommon, resulting in unfortunate situations like the incomparability of $1$ and $0$. 

\item Addition is less intuitive. Disjunctive sum is defined in mis\`ere as in normal play, but much of our intuition for the interaction of games in a sum is lost. For example, the sum of two left-win games may be right-win! In fact, nothing can be said about the addition table of outcomes in misere play:  for any three outcomes $\mathscr{O}_1,\mathscr{O}_2,\mathscr{O}_3 \in \{\mathscr{L, N, P, R}\}$, we can find positions $G$ and $H$ such that $G$ has mis\`ere outcome $\mathscr{O}_1$, $H$ has mis\`ere outcome $\mathscr{O}_2$, and $G+H$ has mis\`ere outcome $\mathscr{O}_3$ \cite{MesdaO2007}.  
Other problems arise with sums, due to the lack of simplification under mis\`ere play: for example, the sum of a game with value $n\in \mathbb{Q}_2$ and a game with value $m\in \mathbb{Q}_2$ may not even be a number-valued game\footnote{This, along with the incomparability of number-valued games, demonstrates that the numerical value system developed for normal play is virtually meaningless in mis\`ere.}, let alone the game with value $n+m$. 

\end{itemize}
For these reasons and others, the study of mis\`ere games was neglected for most of the 20th century. One chapter of {\em On Numbers and Games} presents an analysis of `How to Lose When You Must', and {\em Winning Ways} extends this work in their chapter  `Survival in the Lost World', but both texts consider only impartial mis\`ere games.  The {\em genus theory} developed in the latter allowed for the analysis of certain impartial mis\`ere games, but left most unsolvable \cite{Plamb2008}.  A theory for partizan mis\`ere games seemed, if possible, even more elusive.

The fog began to lift when Thane Plambeck \cite{Plamb2005, PlambS2008} and Aaron Siegel \cite{PlambS2008} introduced a modified equality relation for games under mis\`ere play. Instead of requiring games to be interchangeable in {\em any} sum of games, two games will be considered {\em equivalent modulo} $\mathcal{U}$ if they can be interchanged in any sum of games from the set $\mathcal{U}$.
For example, we might take $\mathcal{U}$ to be the set of all positions that occur in some particular game, such as domineering, and then two domineering positions are equivalent `modulo domineering' if they are interchangeable in any sum of domineering positions. This is a natural and practical  equivalence relation, and its introduction has encouraged renewed interest in the study of mis\`ere games.

Although initially designed only for impartial games, this `restricted' mis\`ere analysis works equally well for partizan games \cite{Siege}. The study of restricted partizan  mis\`ere games began with the doctoral theses of Paul Ottaway \cite{Ottaw2009} and Meghan Allen \cite{Allen2009}, and has continued with a relative flurry of recent activity from a number of additional researchers. The present survey of partizan mis\`ere game theory  highlights the most significant results from recent research, including canonical forms of partizan mis\`ere games, the invertibility of games under restricted mis\`ere play, and applications to specific partizan mis\`ere versions of Nim, Kayles, and Hackenbush.
We begin with some prerequisite definitions.

\section{Prerequisites}

We use the notation $G=\{G^{\lopt} \v G^{\ropt} \}$, where $G^{\lopt}=\{G^{L_1},G^{L_2},\ldots\}$ is the set of left options from $G$ and $G^{L}$ is a particular left option. Any position which can be reached from a game $G$ is called a {\em follower} of $G$.

The {\em outcome} of a game is $\mathscr{L}$ if Left wins playing first or second, $\mathscr{R}$ if Right wins playing first or second, $\mathscr{N}$ if either player can win going first, and $\mathscr{P}$ if neither player can win going first. These outcomes are partially ordered as in normal play: that is, $\mathscr{L}>\mathscr{P}>\mathscr{R}$,  $\mathscr{L}>\mathscr{N}>\mathscr{R}$, and $\mathscr{P} || \mathscr{N}$. We use the outcome function $o^-(G)$ to denote the mis\`ere outcome of $G$ and $o^+(G)$ to denote the normal outcome of $G$. The {\em outcome classes} $\L^-, \N^-, \R^-, \P^-$ are the sets of all games with the indicated outcome under misere play, so that we can write $G\in \L^-$ when $o^-(G)=\mathscr{L}$.

In normal play, the \textit{negative} of a game is defined recursively as $-G=\{-{G}^{\ropt} | 
{-G}^{\lopt}\}$, and is so-called because $G+ (-G)=0$ for all games $G$ under normal play.  As mentioned in the introduction, this property holds in mis\`ere play only if $G$ is  the zero game $\{\cdot \v \cdot\}$  \cite{MesdaO2007}.  To avoid confusion and inappropriate cancellation, we generally write $\overline{G}$ instead of $-G$ and refer to this game as the {\em conjugate} of $G$.
  
Most other definitions from normal-play game theory are used without modification for \mise games, including  disjunctive sum, equality, and inequality. See \cite{AlberNW2007} for an excellent introduction to normal play. In this paper, when equality and inequality relations are used, mis\`ere play is assumed unless otherwise stated.
The equivalence relation developed by Plambeck and Siegel is formalized in Definition~\ref{def:eq} below.

\begin{defn}
\label{def:eq}
For games $G$ and $H$ and a set of games $\mathcal{U}$, the terms 
{\em equivalence} and {\em inequality}, {\em modulo} $\mathcal{U}$, are  defined by
$$G\equiv H \textrm{ (mod } \mathcal{U}) \textrm{ if and only if } o^-(G+X)= o^-(H+X) \textrm{ for all games } X \in \mathcal{U},$$
$$G\geqq H \textrm{ (mod } \mathcal{U}) \textrm{ if and only if } o^-(G+X)\geq o^-(H+X) \textrm{ for all games } X \in \mathcal{U}.$$
\end{defn}

The word {\em indistinguishable} is sometimes used instead of {\em equivalent}, and if $G\not \equiv H$ (mod $\mathcal{U}$) then $G$ and $H$ are said to be {\em distinguishable} modulo $\mathcal{U}$.  In this case there must be a game $X\in \mathcal{U}$ such that $o^-(G+X)\not = o^-(H+X)$, and we say that $X$ {\em distinguishes} $G$ and $H$. 
The set  $\mathcal{U}$ is called the {\em universe}. All universes in this survey are closed under followers and disjunctive sum, and most are also closed under conjugation. Although we usually assume $G$ and $H$ are games in $\mathcal{U}$, this stipulation is unnecessary, and it is sometimes useful to compare games modulo a universe $\mathcal{U}$ even when the games do not belong to $\mathcal{U}$.

Notice that $G\equiv H$ (mod $\mathcal{U}$) implies $G\equiv H$ (mod $\mathcal{V}$) for any subset $\mathcal{V} \subseteq \mathcal{U}$, but in general games can be equivalent in the smaller universe and distinguishable in the larger. Also note that this equivalence is actually a congruence relation with respect to disjunctive sum of games.

\subsection{Specific universes and properties}

A number of  specific game universes are discussed in the sections to follow, and we will define them here. Firstly, we identify games where one or both players have no move:  
a  {\em left end} is a position with no first move for Left (that is, $G$ with $G^{\lopt}=\emptyset$),  a {\em right end} is a position with no first move for Right ($G^{\ropt}=\emptyset$), and an {\em end} is a position that  is either a left end or a right end or both (the zero game). 

A left (right) end is called {\em dead} if each of its followers is also a left (right) end. Games in which every end follower is a dead end are called {\em dead-ending}. Figures 1 and 2 provide examples to illustrate these definitions. By definition, in dead-ending games, if Left has no move at some point, then Left will never have a move again. This is a natural property held by well-studied games such as hackenbush, domineering, and other so-called {\em placement} games (where players move by placing pieces on a board).
 The set of all dead-ending games is denoted $\mathcal{E}$ and has proven to be rich in interesting results for mis\`ere play.  The set of all dead ends and sums of dead ends is denoted $\mathcal{E}_e$.

 \begin{figure} \label{ends}
 \centering
\includegraphics[scale=0.35]{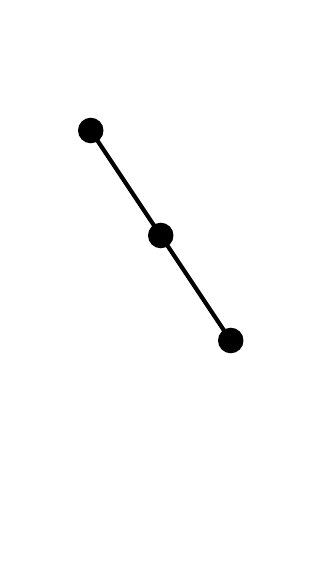} \hspace{1em} \includegraphics[scale=0.35]{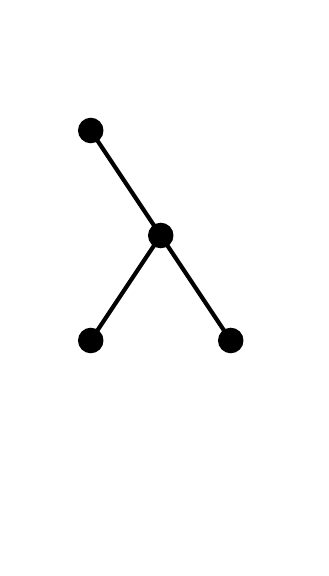}
 \vspace{-1cm}

 \caption{A dead left end and a left end that is not dead.}
 \end{figure}

 \begin{figure} \label{deadending}
 \centering
\includegraphics[scale=0.35]{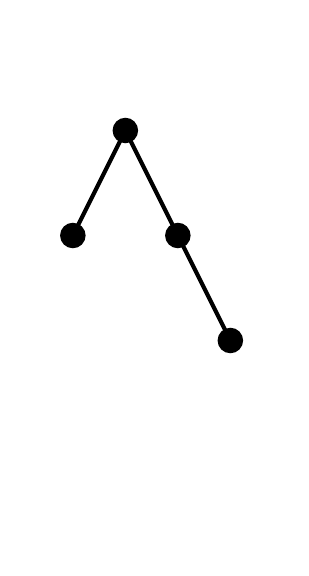} \hspace{1em} \includegraphics[scale=0.35]{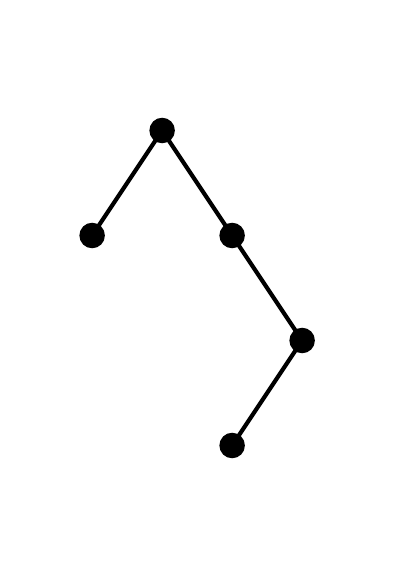}
  \vspace{-1cm}
\caption{A dead-ending game and a game that is not dead-ending.}
\end{figure}
 
 Games in which the only end is zero  --- that is, where Left can move if and only if Right can move --- are called {\em all-small}  in normal play and {\em dicot}  in mis\`ere. The set of all dicot games is denoted $\mathcal{D}$.
  Note that $\mathcal{D}$ is a proper subset of $\mathcal{E}$.

A position is called {\em alternating} if neither player can make consecutive moves; that is, if $G^{L\lopt}$ and $G^{R\ropt}$ are empty for all $G^L$ and all $G^R$. This restriction allows for easier analysis under mis\`ere play. The set of all sums of alternating games is denoted $\mathcal{A}$, and the set of all alternating ends and their sums is denoted $\mathcal{A}_e$. 

We  end this section by mentioning two significant properties that universes may have. The following definitions appear in recent work by Larsson, Nowakowski, and Santos, as part of their new framework called  `absolute game theory', in which they generalize the theories of normal, mis\`ere, and other types of play. Note that their universes are closed under conjugates. We give the definition of {\em density} specifically for the mis\`ere case, but it can be defined generally as well.

\begin{defn} {\em \cite {LarssonNS2016A}}
A universe $\mathcal{U}$ is {\bf parental} if for any two non-empty sets  $\mathcal{A}, \mathcal{B} \subseteq \mathcal{U}$, the position $\{ \mathcal{A} | \mathcal{B} \}$ is also in $\mathcal{U}$.
\end{defn}

\begin{defn} {\em \cite {LarssonNS2016A}}
 A universe $\mathcal{U}$  is {\bf dense} under mis\`ere  play  if, for all $G \in \mathcal{U}$ and any outcome $\mathscr{O}$ in $\{\mathscr{L},\mathscr{R},\mathscr{N},\mathscr{P}\}$, there is an $H\in \mathcal{U}$ such that the mis\`ere outcome of $G+H$  is $\mathscr{O}$.
\end{defn}

These properties are relevant to the current research areas of comparability and invertibility of mis\`ere games. These areas are discussed in Sections 3 and 4; Section 5 discusses the problems with and very recent solutions to the reversibility of mis\`ere games, and Section 6 discusses solutions to specific games under mis\`ere play. We begin with the comparability of mis\`ere games.

\section{Comparability}

In normal play, $G\geq H$ if and only if $G-H$ is previous-win, and so there is an easy test for inequality of games. In general mis\`ere play, we do not have this test, and so proving $G\geq H$ in mis\`ere play requires proving $o(G+X)\geq o(H+X)$ for all games $X$. 

We do at least have a slightly-modified {\em hand-tying} principle for mis\`ere games  \cite{MesdaO2007}.  In normal-play, this principle says that if two games $G$ and $H$ differ only by the addition of one or more extra left options to $G$, then Left can do at least as well playing $G$ as playing $H$ ($G\geq H$, in normal play); at worst, Left can `tie her hand' and ignore the extra options, thereby essentially playing the game $H$ instead of $G$.  In mis\`ere play, the same argument holds, with one stipulation: the set $H^{\lopt}$ of left options cannot be empty.  If it is, adding a left option is not always beneficial to Left, who is sometimes happy to have no move in a position. However, when there already exists at least one left option, Left can simply ignore any additional ones. 
This principle was used in \cite{Milley2013} and \cite{DorbecRSS} to classify day-2 and day-3 dicot games.

 If we restrict ourselves to a particular universe of games $\mathcal{U}$, then we need only consider games $X$ in $\mathcal{U}$, and so  we may be able to show $G\geqq H$ (mod $\mathcal{U}$) even if $G\not \geq H$ in general.
When current research from absolute game theory \cite{LarssonNS2016A} (see also \cite{Ettinger}) is applied to mis\`ere games, we see that comparability of games $G$ and $H$  can be demonstrated  without considering the sum of $G$ and $H$ with {\em all} $X\in \mathcal{U}$, provided certain conditions are met by $\mathcal{U}$. Specifically, if $\mathcal{U}$ is parental and dense, then the following result holds.

\begin{thm} {\em \cite{LarssonNS2016A}} Let $\mathcal{U}$ be a universe that is conjugate-closed, parental, and dense. Then $G\geqq H$ (mod $\mathcal{U}$) if and only if
\begin{enumerate}
\item[(i)]  For all $G^R$ there is $G^{RL} \geqq H$ (mod $\mathcal{U}$) or $G^R \geqq H^R$ (mod $\mathcal{U}$); 
\item[(ii)] For all $H^L$ there is $G^L \geqq  H^L$ (mod $\mathcal{U}$) or $G \geqq H^{LR}$ (mod $\mathcal{U}$);
\item[(iii)] If H is a left end, then Left wins playing first in $G+X$ for any left end $X$ in $\mathcal{U}$; and
\item[(iv)]  If G is a right end, then Right wins playing first in $H+X$ for any right end $X$ in $\mathcal{U}$.
\end{enumerate}

\end{thm}

We will see this importance of ends in other areas of mis\`ere analysis, including invertibility, where we look next.

\section{Invertibility}

As stated in the introduction, no non-zero game has an additive inverse in general mis\`ere play. However, in a restricted universe $\mathcal{U}$, a game $G$ may satisfy $G+\overline{G}\equiv 0$ (mod $\mathcal{U}$) --- or perhaps even $G+H \equiv 0$  (mod $\mathcal{U}$) for $H\not \equiv \overline{G}$  (mod $\mathcal{U}$),  as discussed in Section 4.1 --- and then the game $G$ is said to be {\em invertible} modulo $\mathcal{U}$.
 The first result of this kind was Meghan Allen's demonstration that $*+*\equiv 0$ in any universe of dicot games \cite{Allen2011}. Allen's result is generalized in \cite{MckayMN} with the following  sufficient condition for invertibility in the universe of dicots.

\begin{thm} \label{sprigs} {\em \cite{MckayMN}} If $G+\overline{G} \in \mathcal{N}^-$ and  $H+\overline{H} \in \mathcal{N}^-$ for all followers $H$ of $G$, then $G+\overline{G}\equiv 0 (\textrm{mod } \mathcal{D})$.
 \end{thm}
 
 Theorem~\ref{sprigs} was used to show that the ordinal sum of $*$ and a number\footnote{By {\em number} ({\em integer}) in mis\`ere play, we mean a game that is identical to the normal-play canonical form of a number (integer).}, $*$:$x$, is invertible in the universe of dicots. This  result and others are presented in Table~\ref{tab:inv}, which lists some of the positions known to be invertible in the universes of alternating games ($\mathcal{A}$), dicots ($\mathcal{D}$), or dead-ending games ($\mathcal{E}$).
 
 \begin{table}
\label{tab:inv}
 \centering
 \begin{tabular}[h]{lll}
{\bf Invertible position} & {\bf Universe of invertibility} & {\bf Reference} \\ 
$*$ &  $\mathcal{D}$ &  \cite{Allen2011} \\
$*:x$ for $x\in \mathbb{Q}_2$ & $\mathcal{D}$ & \cite{MckayMN}\\
Any dead end (e.g., $n\in \mathbb{Z}$) & $\mathcal{E}$  & \cite{MilleyRenault}\\
Any alternating end &  $\mathcal{A}$  &\cite{MilleyNO}\\
Any alternating game not in $ \mathcal{P}^-$ &  $\mathcal{A}$ &   \cite{Milley2013}\\

 \end{tabular}
 \caption{Some known instances of invertibility in restricted mis\`ere play.}
 \end{table}

Many of these instances of invertibility were demonstrated using the  following sufficiency condition for invertibility in restricted mis\`ere play. Generally, one proves $G+\overline{G}\equiv 0$ (mod $\mathcal{U}$) by showing that the outcome of $G+\overline{G}$ is the same as the outcome of $G+\overline{G}+X$ for any $X$ in $\mathcal{U}$. Theorem \ref{thesis} essentially says that you need only check the $X$ positions that are ends, an idea that is paralleled by the more recent result about comparability of mis\`ere games (Theorem 3.1).

\begin{thm} {\em \cite{MilleyRenault}}\label{thesis}
Let $U$ be a universe of games closed under followers, sum, and conjugation, and let $S\subseteq U$ be a set of games closed under followers.  
If $G+\overline{G} + X\in \L^-\cup \N^-$ 
for every game $G \in S$ and every left end $X \in \mathcal{U}$, 
 then $G+\Gb \equiv 0$ (mod $\mathcal{U}$) for every $G\in S$.
\end{thm}

\subsection{Non-conjugal invertibility}

A bizarre property of restricted mis\`ere play is that a game $G$ can have an additive inverse modulo some universe $\mathcal{U}$ without that  inverse being the conjugate $\overline{G}$. The only known partizan result of this kind appears in \cite{Milley2015}, where the games $\{0 | \cdot\}$ and $\{1 | 0\}$ sum to zero among the set of all {\em partizan kayles}\footnote{The paper \cite{Milley2015} solves a partizan version of the game Kayles, played on rows of pins, where Left can knock down a single pin and Right can knock down two adjacent pins.} positions, despite neither being equivalent to the conjugate of the other in this universe. This inverse pair is further remarkable for the fact that one position is right-win and the other is previous-win.

In the example from partizan kayles, the actual conjugates of $\{0 | \cdot\}$ and $\{1 | 0\}$ do not even belong to the universe. In \cite{Milley2013} it was conjectured that being closed under conjugation would prevent such occurrences of non-conjugal invertibility; however, a counterexample can be seen in Appendix 6 of \cite{PlambS2008}, for a subset of impartial games.

This leads us to a pressing open question in mis\`ere theory: in what universes $\mathcal{U}$ do we have no non-conjugate inverses, so that $G+H\equiv 0$ (mod $\mathcal{U}$) only if $H\equiv \overline{G}$ (mod $\mathcal{U}$)? In-progress research suggests that adapting the proof of a similar result on scoring games \cite{LarssonNPS} can prove that no non-conjugate inverses occur in  dicot games,  dead-ending games, or any universe that is parental, dense, and amenable to a  type of `replacement' reversibility through ends, as discussed in the next section.

\section{Reversibility and canonical forms}

Given the relative lack of structure in mis\`ere play, it is perhaps surprising that we have canonical forms here just as in normal play, with precisely the same definitions of domination and reversibility (with inequality under mis\`ere play instead of normal play). This was shown in the collaborative paper of G.A. Mesdal \cite{MesdaO2007} and subsequent work by Aaron Siegel \cite{Siege}. The latter  also demonstrated that, as in normal play, the simplified game obtained by removing dominated options and bypassing reversible ones is unique.

So canonical forms `work' in mis\`ere play; but in general the concept is less useful than in normal play, because it is so hard to find instances of domination or reversibility. If we restrict ourselves to a particular universe of games $\mathcal{U}$, then we may get domination or reversibility in the restricted universe that does not occur in general, due to inequalities of the form  $G\geqq H$ (mod $\mathcal{U}$). Consequently, a game could have different `restricted canonical forms' in different universes. However, the construction of  a canonical form --- specifically, how we deal with reversible options --- is not quite the same when the universe is restricted in this way.

 The problem of reversibility  in restricted universes is related to the following result of \cite{Siege}, which is used in the construction of mis\`ere canonical forms.

\begin{lem} \em{[\cite{Siege}, Lemma 3.5]}
If $H$ is a left end and $G$ is not, then $G\not \geq H$.
\end{lem}

This result holds in the context of all mis\`ere games; however, it may be the case that a non-Left-end $G$ can be greater than a left end $H$ {\em modulo some universe} $\mathcal{U}$. For example, in the universe of dicot games $\mathcal{D}$, we have $ \{0,* | *\} \geq 0$ (mod $\mathcal{D}$) \cite{DorbecRSS}.

Why is this a problem? In general, Lemma 3.5 means we never have to worry about reversibility through an end; it cannot happen that $G\geq G^{LR}$ if $G^{LR}$ is a left end, and so in such a case $G^L$ could not be reversible. This fact is exploited in the proof that reversibility works in mis\`ere play: that $G'=G$ when $G'$ is obtained from $G$ by replacing a reversible option $G^L$ with the left options of $G^{LR}$ \cite{Siege}. Since the same fact does not hold in restricted mis\`ere play, the result from \cite{Siege} no longer applies, and so we cannot necessarily bypass all reversible options. In the example above from \cite{DorbecRSS}, even though $G=\{0,* | *\} \geq 0$ (mod $\mathcal{D}$) and $0= *^R = G^{LR}$, it is not the case that $\{0,* | *\} \equiv \{0|*\}$ (mod $\mathcal{D}$). Left's only good move in $G$ is to $*$, so removing $*$ and replacing it with no options does not result in a game that is just as good for Left. 

In \cite{DorbecRSS}, the proof from \cite{Siege} of uniqueness of mis\`ere canonical forms was adapted to construct unique restricted canonical forms in the universe of dicot games. For dicots, the problem of reversibility through ends is dealt with as follows:  if $G^L$ is reversible through a left end, then replacing $G^L$ with $*$ results in an equivalent game. This solution should be further adaptable to other restricted universes: we would just need to find a suitable `replacing game', that might depend on $G$, to replace an option $G^L$ that is reversible  through a left end. As the invertibility of $*$ (modulo dicots) is used in the proof of the uniqueness of the  canonical forms for dicots, the replacing game in other universes will most likely  have to be invertible. Solving the problem of reversibility in specific mis\`ere universes is an open area of research; notably, in-progress work from the authors of \cite{LarssonNS2016A} suggests a solution for certain universes, including dead-ending games.

This completes our survey of the recent developments in mis\`ere theory, including comparability, invertibility, and reversibility of mis\`ere games. We next show how some of these advancements have been applied to solve specific games under mis\`ere play.

\section{Applications to specific games}
 
A number of specific partizan games have been successfully solved using the theory of restricted mis\`ere play. These solutions consider equivalence modulo the universe of all positions that occur under the specific game rule set, and take advantage of results for broader superset universes.

Penny Nim is a partizan variant of Nim played with stacks of coins. 
In each stack, coins are all heads up or all tails up, and the entire stack may be lying sideways.
On her turn, Left chooses a stack with tails-up or sideways orientation, removes any number of coins from it, and turn it heads up.
Right plays similarly on heads-up or sideways coins stacks, but leaves them tails-up. 
Notice that any position of this game is alternating, and the potential for sideways stacks means that not all components are initially ends.
The game is solved in \cite{Milley2013}, using the analysis of the alternating universe $\mathcal{A}$, in which most `single-stack' positions are invertible. The solution involves first simplifying single stacks of coins,  modulo $\mathcal{A}$, and then determining outcomes of sums of these finitely-many simplified positions.

Partizan Kayles is a variant of Kayles, played on a row of pins, where Left can knock down a single pin and Right can knock down exactly two adjacent pins.
Notice that any position of this game is dead-ending.
The game is solved in~\cite{Milley2015}. The key is to see that Left should always take an isolated single pin when she can; this allows for removal of dominated options and decomposition of long rows of pins into shorter rows --- into only isolated single pins and pairs of pins, in fact --- and then all that remains is to see who wins on a sum of such positions. This is easily done once it is shown that an isolated single pin and an isolated pair of adjacent pins `cancel' (that is, they are additive inverses). 

Hackenbush Sprigs is a particular case of the game of Hackenbush.
The game can be seen as rows of blue, green and red dominoes where each row has exactly one green domino, which is the leftmost domino.
A move of Left is to pick a blue or green domino and remove it with all dominoes of the same row to its right.
Right plays similarly with red or green dominoes.
Notice that any position of this game is dicot.
The game is solved in~\cite{MckayMN}.
The authors first show that all games are invertible by finding the canonical forms of all rows, modulo dicot games, and then finding the outcomes of sums of such positions.
They end by showing no other simplification can be made, completely solving the game.

\section{Current and future directions}

We conclude by highlighting two of the open problems that were introduced above. 
\begin{enumerate}
\item {\bf Non-conjugal invertibility.} In what universes does $G+H\equiv 0$ imply $\overline{G} \equiv H$?  Can we indeed prove that this is true for certain parental, dense universes, and if so, what known games naturally occur in such universes? Can we find other examples of universes in which this is {\em not} true, besides the one impartial and one partizan example that have been identified to date?

\item {\bf Reversibility through ends.} There is a solution for dicots, where options that are reversible through ends are replaced with $*$, and there is a suggestion that a similar solution will work in a few other specific universes. Can we solve reversibility in  other universes, perhaps starting even with small ruleset-specific universes? Can we find a general process for constructing the necessary `replacement' games? Are there universes in which the problem of reversibility through ends does not even arise --- that is, in which Lemma 5.1 holds?
\end{enumerate}

\end{document}